\documentclass[amsbook,12pt,a4paper]{article}
\usepackage{epsfig, amsfonts, amssymb, amsmath}
\usepackage[usenames]{color}
\usepackage{picinpar}
\newcommand{\nl}{\mbox{}\\}

\newcommand{\bb}{\mbox{\boldmath $b$}}
\newcommand{\uu}{\mbox{\boldmath $u$}}
\newcommand{\RR}{\mathbb{R}}
\newcommand{\LL}{\mbox{\boldmath $L$}}
\newcommand{\HH}{\mbox{\boldmath $H$}}
\begin{document}
\setcounter{page}{1}
%
% ----------------------------------------------------------------
%
\mbox{} \vspace{-2.000cm} \\
\begin{center}
{\Large \bf
 Decay Rates of global weak solutions} \\
\mbox{} \vspace{-0.300cm} \\
{\Large \bf
for the MHD equations in  $\dot{\HH^s}(\RR^n)$} \\
% in \mbox{\boldmath $\mathbb{R}^{n}\!\:\!$}} \\
%
%
% ------------------------------------------------
%
\nl
\mbox{} \vspace{-0.300cm} \\
{\sc Robert Guterres, Juliana Nunes and Cilon Perusato } \\
\mbox{} \vspace{-0.350cm} \\
{\small Departamento de Matem\'atica Pura e Aplicada} \\
\mbox{} \vspace{-0.670cm} \\
{\small Universidade Federal do Rio Grande do Sul} \\
\mbox{} \vspace{-0.670cm} \\
{\small Porto Alegre, RS 91509, Brazil} \\
\nl
\mbox{} \vspace{-0.300cm} \\
%
%
% ------------------------------------------------
%
%
{\bf Abstract} \\
\mbox{} \vspace{-0.300cm} \\
\begin{minipage}{12.70cm}
{\small
\mbox{} \hspace{+0.250cm}
We show that
$ t^{s/2} \| (\uu,\bb)(\cdot,t) \|_{{\dot{H}^s(\RR^n)}} \to 0 $ as $t \to \infty$
for Leray solutions $(\uu,\bb)(\cdot,t)$
of the incompressible MHD equations,
where $ 2 \leq n \leq 4$ and $ s \geq 0$.
As a corollary of main result described previously
we have also that
$\lim_{t \to \infty} t^{\frac{n}{4} - \frac{n}{2q}}\| (\uu,\bb)(\cdot,t) \|_{\mathbf{L}^q(\RR^n)} = 0$, $2 \leq q \leq \infty$.
}
\end{minipage}
\end{center}

\nl
\mbox{} \vspace{-0.550cm} 
{\small 35D30},  \\
{\small 76W05} \\
\nl
\mbox{} \vspace{-0.800cm} \\
{\bf Key words}: decay rates in Sobolev homogeneous spaces, MHD incompressible, time decay of solution derivatives, Leray global weak solutions, $L^p$ estimates
%incompressible Navier-Stokes equations,
%Leray's (weak) solutions, \linebreak
%\mbox{} \hfill
%large time behavior,
%Leray's $L^{2}\!\;\!$ conjecture,
%energy estimates, heat semigroup. \\
%
%
\\
%
% ----------------------------------------------------------
%

%
% ************************************************
% *                                              *
% *                Section 1                     *
% *              Introduction                    *
% *                                              *
% ************************************************
%
\nl
\mbox{} \vspace{-0.650cm} \\

{\bf 1. Introduction} \\
\setcounter{section}{1}
\par In this work we derive a general decay rate for Leray global  weak solutions of incompressible MHD equations (in $ \dot{\HH^s}(\RR^n)$, where $ n\leq 4$), that is,
global solutions
$(\uu,\bb)(\cdot,t) \in L^{\infty}((0,\infty), \:\!\mbox{\boldmath $L$}^{2}_{\sigma}(\mathbb{R}^{n}))
\,\cap
 $
$ {\displaystyle
	L^{2}((0,\infty), \:\!\dot{\mbox{\boldmath $H$}}\mbox{}^{\!\;\!1}\!\;\!(\mathbb{R}^{n}))
	\cap \;\!
	C_{w}([\;\!0, \infty), \mbox{\boldmath $L$}^{2}(\mathbb{R}^{n}))
} $ of the system
\begin{equation*}
\tag{1.1$a$}
\begin{split}
\mbox{\boldmath $u$}_t
\;\!+\,
\mbox{\boldmath $u$} \cdot \nabla \mbox{\boldmath $u$}
\,+\;\!
\nabla \:\!{P}
\;=\;
\mu \, \Delta \mbox{\boldmath $u$}
\,+\,
\mbox{\boldmath $b$} \cdot \nabla \mbox{\boldmath $b$},
\end{split}
\end{equation*}
\begin{equation*}
\tag{1.1$b$}
\begin{split}
\mbox{\boldmath $b$}_t
\;\!+\,
\mbox{\boldmath $u$} \cdot \nabla \mbox{\boldmath $b$}
\;=\;
\nu \, \Delta \:\!\mbox{\boldmath $b$}
\,+\,
\mbox{\boldmath $b$} \cdot \nabla \mbox{\boldmath $u$},
\end{split}
\end{equation*}
\begin{equation*}
\tag{1.1$c$}
\begin{split}
\nabla \cdot \mbox{\boldmath $u$}(\cdot,t) \,=\, 0,
\quad \;\;\,
\nabla \cdot \mbox{\boldmath $b$}(\cdot,t) \,=\, 0,
\end{split}
\end{equation*}
\mbox{} \vspace{-0.200cm} \\
with initial data
$ {\displaystyle
	(\:\!\mbox{\boldmath $u$}_0, \mbox{\boldmath $b$}_0) \in
	\mbox{\boldmath $L$}^{2}_{\sigma}(\mathbb{R}^{n})
	\!\times\!
	\mbox{\boldmath $L$}^{2}_{\sigma}(\mathbb{R}^{n})
} $,
that $\| (\uu,\bb)(\cdot,t) - (\uu_0,\bb_0) \|_{L^2(\RR^n)}  \to 0$
as $t \to 0$
and such that the \textit{strong} energy inequality\footnote{%
For the definition of $\|(\uu,\bb) \|_{L^2(\RR^n)} $ , $\|(D\uu,D\bb) \|_{L^2(\RR^n)} $
and other similar expressions throughout the text, see (1.7$e$) and (1.7$f$).
}
\begin{equation*}
\tag{1.2}
\begin{split}
\| (\uu,\bb)(\cdot,t) \|^2_{L^2(\RR^n)} + 2\mu \int_{r}^{t} {\| D\uu (\cdot,\tau) \|}^2_{L^2(\RR^n)} d\tau + 2\nu \int_{r}^{t} {\| D\bb (\cdot,\tau) \|}^2_{L^2(\RR^n)} d\tau \\
\leq\;
\| (\uu,\bb)(\cdot,r) \|^2_{L^2(\RR^n)},  \text{   }\forall t \geq r   &
\end{split}
\end{equation*}
holds for a.e $ r \geq 0$, including $ r = 0$.
Such solutions were first constructed by Leray (\cite{Leray1934})
for the Navier-Stokes system where $n \leq 3$
and later by other authors with different methods considering also $n = 4$
or even in higher dimensions,
see e.g.\;\cite{FujitaKato1964,Kato1984,Ladyzhenskaya1969,Lions1969}.
All these methods can be adapted for the MHD equations
\cite{G.DUVAUT&J.L.LIONS,SchonbekSchonbekSuli1996,SermangeTemam1983}.
In (1.1),
$ \mu, \nu > 0 $ are given constants,
\mbox{$ \mbox{\boldmath $u$} = \mbox{\boldmath $u$}(x,t) $}, \mbox{$ \mbox{\boldmath $b$} = \mbox{\boldmath $b$}(x,t) $}
 and $ P = P(x,t) $
are the unknowns
(the flow velocity, magnetic field and total pressure, respectively).
As usual,
$ \!\;\!\mbox{\boldmath $L$}^{2}_{\sigma}(\mathbb{R}^{n}) $
is the
space
of solenoidal fields
$ \:\!\mbox{\bf v} = (v_{1}, ..., v_{n}) \!\:\!\in
\mbox{\boldmath $L$}^{2}(\mathbb{R}^{n}) \equiv L^{2}(\mathbb{R}^{n})^{n} \!\:\!$
with
$ \nabla \!\cdot \mbox{\bf v} \!\;\!= 0 $
% \mbox{in $ {\cal D}^{\prime}(\mathbb{R}^{n}) $}
% (i.e., as distributions)
in the distributional sense,
$ \dot{\mbox{\boldmath $H$}}\mbox{}^{s}(\mathbb{R}^{n}) =
\dot{H}^{s}(\mathbb{R}^{n})^{n} $
where
$ \dot{H}^{s}(\mathbb{R}^{n}) $
denotes the homogeneous Sobolev space of order~$s \geq 0$,
and
$ C_{w}(I, \:\!\mbox{\boldmath $L$}^{2}(\mathbb{R}^{n})) $
denotes the set of mappings from a
given interval~\mbox{$ I \subseteq \mathbb{R} $}
to
$ \mbox{\boldmath $L$}^{2}(\mathbb{R}^{n}) $
that are $L^{2}$-\;\!weakly continuous
at each $\:\! t \in I$.
Here,
we always assume $ 2 \leq n \leq 4 $.
 Moreover, similarly to the Navier-Stokes case,
 there always exists\footnote{If $n=2$, then $t_* = 0$.} 
 some $ t_* \gg 1$ 
 -- depending on the solution $(\uu,\bb) $ -- such that one has
 \begin{equation*}
 \tag{1.3$a$}
(\uu,\bb) \in C^\infty(\RR^n \times [t_*,\infty))
 \end{equation*} and, for each $ m \in \mathbb{Z}_+$:
\begin{equation*}
\tag{1.3$b$}
(\uu,\bb)(\cdot,t) \in \LL^\infty ([t_*,T), \HH^m(\RR^n)) ,
\end{equation*}
for each $t_* < T < \infty$, that is, $ (\uu,\bb)(\cdot,t) \in \LL^\infty_{\text{loc}} ([t_*,\infty), \HH^m(\RR^n))$.
%
%----------------------------------------------------------
%
%         Introduction: Beginning
%
%----------------------------------------------------------
%
In \cite{Schonbek&RubenAgapito2007}, 
Agapito and Schonbek showed that
\begin{equation*}
\tag{1.4}
\|(\uu,\bb)(\cdot,t)\|_{L^2(\RR^n)} \to 0 \text{   as   } t \to \infty,
\end{equation*}
generalizing the Kato techniques for Navier-Stokes equation (see e.g \cite{Kato1984}) 
in dimension $n=2,3$. 
More recently, in \cite{SchutzZinganoZingano2015} 
the authors showed the above property for the Navier-Stokes equations 
in a simple way using Duhamel's principle and, with the same technique, 
they provided an $L^\infty$ decay rate. 
So, we will adapt this for the MHD equations in a preliminaries section 
and generalize this argument to obtain (1.4) in $n=4$ dimension. 
However, it was necessary to prove the following decay property for derivatives
\begin{equation*}
\lim_{t\,\rightarrow\,\infty}
\,
t^{\:\!1/2}
\,
\|\, (D\uu,D\bb)(\cdot,t) \,
\|_{\mbox{}_{\scriptstyle L^{2}(\mathbb{R}^{n})}}
=\; 0,
\quad \;\;\,
2 \leq n \leq 4.
\end{equation*}
%Starting\footnote{For $n=2$, the solution of (1.1) is always smooth, in other words, $t_* = 0$}
%with $n = 3$: we choose a $t_0 \geq t_*$ (see 1.3).
% Since $(\uu,\bb)(\cdot,t)$ is smooth for large t, applying the Duhamel's principle to the system (1.1), one has
Studying these problems with this new approach 
we provide a decay rate for all the derivatives 
and using interpolation 
we get the general decay below.
\\
\\
\textbf{Main Theorem.}
\\
\textit{For a Leray solution $(\uu,\bb)(\cdot,t)$ of (1.1) and $ n \leq 4$, one has}
\begin{equation*}
\tag{1.5}
\begin{split}
\lim_{t \to \infty} t^{s/2} \| (\uu,\bb)(\cdot,t) \|_{{\dot{H}^s(\RR^n)}} = 0,
\end{split}
\end{equation*}
for all $s \geq 0$. \\
%}
%
\nl
As a consequence, we get the following result. \\
\nl
\textbf{Corollary.}
\\
\textit{For a Leray solution $(\uu,\bb)(\cdot,t)$  of (1.1) and $ n \leq 4$, one has}
\begin{equation*}
\tag{1.6}
\begin{split}
\lim_{t \to \infty} t^{\frac{n}{4} - \frac{n}{2q}}\| (\uu,\bb)(\cdot,t) \|_{{L}^q(\RR^n)} = 0,
\end{split}
\end{equation*}
$2 \leq q \leq \infty$. \\
In Section 2, we recall some basic facts and estimates that are needed (or have relevance) for the derivation of Main Theorem in Section 3. \\
\nl
%
% -------------------------------------
%
%            Notation
%
% -------------------------------------
%
{\bf Notation.}
As shown above,
boldface letters
are used for
vector quantities,
as in
$ {\displaystyle
	\;\!
	\mbox{\boldmath $u$}(x,t)
	=
} $
$ {\displaystyle
	(\:\! u_{\mbox{}_{\!\:\!1}}\!\;\!(x,t), ...\,\!\,\!,
	\:\! u_{\mbox{}_{\scriptstyle \!\;\! n}}\!\;\!(x,t) \:\!)
} $.
$\!$Also,
$ \nabla P \;\!\equiv \nabla P(\cdot,t) $
denotes the spatial gradient of $ \;\!P(\cdot,t) $,
$ D_{\!\;\!j} \!\;\!=\:\! \partial / \partial x_{\!\;\!j} $,
$ {\displaystyle
	\,\!
	\nabla \!\cdot \mbox{\boldmath $u$}
	\:\!=
	D_{\mbox{}_{\!\:\!1}} u_{\mbox{}_{\!\:\!1}} \!\;\!+
	... \!\;\!+
	D_{\mbox{}_{\scriptstyle \!\;\!n}} \,\!
	u_{\mbox{}_{\scriptstyle \!\;\!n}}
} $
is the (spatial) divergence of
$ \:\!\mbox{\boldmath $u$}(\cdot,t) $.
$ |\,\!\cdot\,\!|_{\mbox{}_{2}} \!\,\!\,\!$
denotes the Euclidean norm
in $ \mathbb{R}^{n} \!$,
and
$ {\displaystyle
	\,\!\,\!
	\| \:\!\cdot\:\!
	\|_{\scriptstyle L^{q}(\mathbb{R}^{n})}
	\!\;\!
} $,
$ 1 \leq q \leq \infty $,
are the standard norms
of the Lebesgue spaces
$ L^{q}(\mathbb{R}^{n}) $,
with the vector counterparts \\
\mbox{} \vspace{-0.625cm} \\
\begin{equation}
\tag{1.7$a$}
\|\, \mbox{\boldmath $u$}(\cdot,t) \,
\|_{\mbox{}_{\scriptstyle L^{q}(\mathbb{R}^{n})}}
\;\!=\;
\Bigl\{\,
\sum_{i\,=\,1}^{n} \int_{\mathbb{R}^{n}} \!
|\:u_{i}(x,t)\,|^{q} \,dx
\,\Bigr\}^{\!\!\:\!1/q}
\end{equation}
\mbox{} \vspace{-0.700cm} \\
\begin{equation}
\tag{1.7$b$}
\|\, D \mbox{\boldmath $u$}(\cdot,t) \,
\|_{\mbox{}_{\scriptstyle L^{q}(\mathbb{R}^{n})}}
\;\!=\;
\Bigl\{\,
\sum_{i, \,j \,=\,1}^{n} \int_{\mathbb{R}^{n}} \!
|\, D_{\!\;\!j} \;\!u_{i}(x,t)\,|^{q} \,dx
\,\Bigr\}^{\!\!\:\!1/q}
\end{equation}
\mbox{} \vspace{-0.350cm} \\
and, in general, \\
\mbox{} \vspace{-0.750cm} \\
\begin{equation}
\tag{1.7$c$}
\|\, D^{m} \mbox{\boldmath $u$}(\cdot,t) \,
\|_{\mbox{}_{\scriptstyle L^{q}(\mathbb{R}^{n})}}
\;\!=\;
\Bigl\{\!\!
\sum_{\mbox{} \;\;i, \,j_{\mbox{}_{1}} \!,..., \,j_{\mbox{}_{m}} =\,1}^{n}
\!\;\! \int_{\mathbb{R}^{n}} \!
|\, D_{\!\;\!j_{\mbox{}_{1}}}
\!\!\!\;\!\cdot\!\,\!\cdot\!\,\!\cdot \!\:\!
D_{\!\;\!j_{\mbox{}_{m}}}
\!\:\! u_{i}(x,t)\,|^{q} \,dx
\,\Bigr\}^{\!\!\:\!1/q}
\end{equation}
\mbox{} \vspace{-0.175cm} \\
if $ 1 \leq q < \infty $\/;
if $\, q = \infty $,
then
$ {\displaystyle
	\;\!
	\|\, \mbox{\boldmath $u$}(\cdot,t) \,
	\|_{\mbox{}_{\scriptstyle L^{\infty}(\mathbb{R}^{n})}}
	\!=\;\!
	\max \, \bigl\{\,
	\|\,u_{i}(\cdot,t)\,
	\|_{\mbox{}_{\scriptstyle L^{\infty}(\mathbb{R}^{n})}}
	\!\!: \, 1 \leq i \leq n
	\,\bigr\}
} $, \linebreak
\mbox{} \vspace{-0.530cm} \\
$ {\displaystyle
	\|\, D \,\!\mbox{\boldmath $u$}(\cdot,t) \,
	\|_{\mbox{}_{\scriptstyle L^{\infty}(\mathbb{R}^{n})}}
	\!=\;\!
	\max \, \bigl\{\,
	\|\, D_{\!\;\!j} \;\! u_{i}(\cdot,t)\,
	\|_{\mbox{}_{\scriptstyle L^{\infty}(\mathbb{R}^{n})}}
	\!\!: \:\! 1 \leq i, \:\!j \leq n
	\,\bigr\}
} $
and,
for general \mbox{$m \!\;\!\geq\!\:\! 1$\/:} \\
\mbox{} \vspace{-0.550cm} \\
\begin{equation}
\tag{1.7$d$}
\|\, D^{m} \mbox{\boldmath $u$}(\cdot,t) \,
\|_{\mbox{}_{\scriptstyle L^{\infty}(\mathbb{R}^{n})}}
\!\;\!=\;
\max\,\Bigl\{\;\!
\|\;\! D_{\!\;\!j_{\mbox{}_{1}}}
\!\!\!\;\!\;\!\,\!\cdot \!\;\!\cdot \!\;\!\cdot \!\;\!\;\!
D_{\!\;\!j_{\mbox{}_{m}}}
\!\!\;\!\;\! u_{i}(\cdot,t)\,
\|_{\mbox{}_{\scriptstyle L^{\infty}(\mathbb{R}^{n})}}
\!\!\!\;\!: \;\!
1 \leq \!\;\!\;\!
i, \!\;\!\;\!j_{\mbox{}_{1}}\!\!\;\!\;\!, \!...\!\;\!\;\!,j_{\mbox{}_{m}}
\!\leq n
\!\;\!\;\!\Bigr\}.
\end{equation}
\mbox{} \vspace{-0.200cm} \\
Definitions (1.4) are convenient,
but not essential.
However,
some choice for the vector norms
has to be made to fix the values
of constants. 
We define also, for simplicity the following norms 
for a pair $(\uu,\bb)$ as usually made in literature:
\begin{equation}
\tag{1.7$e$}
\| (\uu,\bb) \|^q_{L^q(\RR^n)} := 
\|\uu\|^q_{L^q(\RR^n)} + \|\bb\|^q_{L^q(\RR^n)}
\end{equation}
and more generally, for all $m \geq 1$ integer
\begin{equation}
\tag{1.7$f$}
\| (D^m\uu,D^m\bb) \|^q_{L^q(\RR^n)} := 
\|D^m\uu\|^q_{L^q(\RR^n)} + \|D^m\bb\|^q_{L^q(\RR^n)}
\end{equation}
for all $1 \leq q \leq \infty$. Similarly, for all $s \geq 0$,
\begin{equation*}
\tag{1.7$g$}
\|(\uu,\bb)\|^2_{\dot{H}^s(\mathbb{R}^{n})} := \|\uu\|^2_{\dot{H}^s(\mathbb{R}^{n})} +
\|\bb\|^2_{\dot{H}^s(\mathbb{R}^{n})},
\end{equation*}
where,
\begin{equation*}
\tag{1.7$h$}
\|\uu\|_{\dot{H}^s(\mathbb{R}^{n})} = \Bigg( \sum_{i=1}^{n} \int_{\RR^n}      |\xi|^{2s}|\hat{u_i}(\xi)|^2 d\xi\Bigg)^{1/2}
\end{equation*}
and $\hat{u_i}$
denote the Fourier transform of $u_i$.
The constants will be represent ed by the letters C, c or K.
For economy, we will typically use the same symbol to denote constants
with different numerical values. \\
%
% ---------------------------------------------------------
%           END OF NOTATION
%%%%%%%%%%%%%%%%%%%%%%%%%%%%%%%%%%%%%%%%%%%%%%%%%%%%%%%%%%%
%
\mbox{} \vspace{-0.300cm} \\
%
%
% ***************************************************
% *                                                 *
% *                 Section 2                       *
% *                                                 %	
% *                Preliminaries                    *
% ***************************************************
%{}
\mbox{} \\
{\bf 2. Preliminaries}
\setcounter{section}{2}
%AQUI começa o comentário
\iffalse
\begin{equation*}
\tag{1.6$a$}
\uu(x,t) = e^{\mu{\Delta(t-t_0)}}\uu(\cdot,t_0) + \int_{t_0}^{t}
e^{\mu\Delta(t-\tau)}  Q_1(\cdot,\tau) d\tau
\end{equation*}
and
\begin{equation*}
\tag{1.6$b$}
\bb(x,t) = e^{\nu{\Delta(t-t_0)}}\bb(\cdot,t_0) + \int_{t_0}^{t}
e^{\nu\Delta(t-\tau)}  Q_2(\cdot,\tau) d\tau,
\end{equation*}
for all $t >t_0 $, where
\begin{equation*}
\tag{1.7$a$}
Q_1 = -\uu\cdot\nabla\uu - \nabla P + \bb\cdot\nabla\bb,
\end{equation*}
\begin{equation*}
\tag{1.7$b$}
Q_2 = -\uu\cdot\nabla\bb  + \bb\cdot\nabla\uu
\end{equation*}
and $e^{\Delta s }$ denotes the heat semigroup.
\fi
%AQUI TERMINA o comentário
\\ \par First, we will obtain the derivatives monotonicity in $L^2(\RR^n)$,
\begin{equation*}
\|(D\uu,D\bb)(\cdot,t)\|_{L^2(\RR^n)} \leq \|(D\uu,D\bb)(\cdot,t_0)\|_{L^2(\RR^n)}.
\end{equation*} Staring with $n = 3$.
This next argument is
adapted from \cite{KreissHagstromLorenzZingano2003}. 
Using  (1.1) and  (1.3), we get,
\begin{multline*}
\tag{2.1}
\|\,  \,\!(D\uu,D\bb)(\cdot,t) \,
\|_{\mbox{}_{\scriptstyle L^{2}(\mathbb{R}^{3})}}
^{\:\!2}
\!\;\!  +\:
2\, \min\{\mu,\nu\} \!\!\;\!
\int_{\scriptstyle t_0}^{\;\!t}
\!\!\:\!
\,
\|\,  \,\!(D^{2}\uu,D^{2}\bb) (\cdot,\tau) \,
\|_{\mbox{}_{\scriptstyle L^{2}(\mathbb{R}^{3})}}
^{\:\!2}
d\tau
\\
\leq\;
\|\,  \,\!(D\uu,D\bb)(\cdot,t_0) \,
\|_{\mbox{}_{\scriptstyle L^{2}(\mathbb{R}^{3})}}
^{\:\!2}
\\ +\;\:\!
C
\!\!\,\!
\int_{\scriptstyle t_0}^{\;\!t}
\|\, (\uu,\bb)(\cdot,\tau) \,
\|_{\mbox{}_{\scriptstyle L^{\infty}(\mathbb{R}^{3})}}
\|\,  \,\!(D\uu,D\bb)(\cdot,\tau) \,
\|_{\mbox{}_{\scriptstyle L^{2}(\mathbb{R}^{3})}}
\|\,  \,\!(D^{2}\uu,D^{2}\bb)(\cdot,\tau) \,
\|_{\mbox{}_{\scriptstyle L^{2}(\mathbb{R}^{3})}}
d\tau
\\
\leq\;
\|\,  \,\!(D\uu,D\bb)(\cdot,t_0) \,
\|_{\mbox{}_{\scriptstyle L^{2}(\mathbb{R}^{3})}}
^{\:\!2}
\\ +\;\:\!
C
\!\!\,\!
\int_{\scriptstyle t_0}^{\;\!t}
\!\!\;\!
\|\, (\uu,\bb)(\cdot,\tau) \,
\|^{1/2}_{\mbox{}_{\scriptstyle L^{2}(\mathbb{R}^{3})}}
\|\,  \,\!(D\uu,D\bb)(\cdot,\tau) \,
\|^{1/2}_{\mbox{}_{\scriptstyle L^{2}(\mathbb{R}^{3})}}
\|\,  \,\!(D^{2}\uu,D^{2}\bb)(\cdot,\tau) \,
\|_{\mbox{}_{\scriptstyle L^{2}(\mathbb{R}^{3})}}
d\tau,
\end{multline*}
where we have used the Sobolev-Nirenberg-Gagliardo (SNG) inequalties (see (2.10$a$)). By (1.2), we can choose $ t_0 \geq t_*$ large enough such that
\begin{equation*}
C^2 \|(\uu_0,\bb_0)\|_{L^2(\RR^3)}
\|(D\uu,D\bb)(\cdot,t_0)\|_{L^2(\RR^3)} < (\min \{\mu,\nu \})^2,
\end{equation*}
so that (2.1) gives $ \|(D\uu,D\bb)(\cdot,t)\|_{L^2(\RR^3)} \leq \|(D\uu,D\bb)(\cdot,t_0)\|_{L^2(\RR^3)} $ for all $t$ near $t_0$ by continuity. Actually, with this choice, it follows from ((2.1) again) that
\begin{equation*}
\tag{2.2}
C^2 \|(\uu_0,\bb_0)\|_{L^2(\RR^3)}
\|(D\uu,D\bb)(\cdot,s)\|_{L^2(\RR^3)} < (\min \{\mu,\nu \})^2, \text{       } \forall s \geq t_0.
\end{equation*}
Recalling (2.1), (2.2) implies that
\begin{equation*}
\tag{2.3}
\|(D\uu,D\bb)(\cdot,t)\|_{L^2(\RR^3)} \leq \|(D\uu,D\bb)(\cdot,t_0)\|_{L^2(\RR^3)},
\end{equation*}
for all $ t \geq t_0$. Since, by (1.2), $\| (D\uu,D\bb)(\cdot,t)\|^2_{L^2(\RR^3)} $ is integrable in $(0,\infty)$ one has, by (2.3),  that\footnote{
	Because a monotonic function $ f \in C^0 ((a,\infty)) \cap L^1((a,\infty)) $ has to satisfy $f(t) = o(1/t)$ as $t \to \infty$ (see e.g. \cite{KreissHagstromLorenzZingano2003}, p. 236).  	}
\begin{equation*}
\tag{2.4}
\lim_{t \to \infty} t \| (D\uu,D\bb)(\cdot,t)\|^2_{L^2(\RR^3)} = 0.
\end{equation*}
A similar argument hold for $ n = 2$ (with $t_* = 0$). For $n = 4$, we proceed as before,
\begin{multline*}
\|\,  \,\!(D\uu,D\bb)(\cdot,t) \,
\|_{\mbox{}_{\scriptstyle L^{2}(\mathbb{R}^{4})}}
^{\:\!2}
\!\;\!  +\:
2\, \min\{\mu,\nu\} \!\!\;\!
\int_{\scriptstyle t_0}^{\;\!t}
\!\!\:\!
\,
\|\,  \,\!(D^{2}\uu,D^{2}\bb) (\cdot,\tau) \,
\|_{\mbox{}_{\scriptstyle L^{2}(\mathbb{R}^{4})}}
^{\:\!2}
d\tau
\\
\leq\;
\|\,  \,\!(D\uu,D\bb)(\cdot,t_0) \,
\|_{\mbox{}_{\scriptstyle L^{2}(\mathbb{R}^{4})}}
^{\:\!2}
\\ +\;\:\!
C
\!\!\,\!
\int_{\scriptstyle t_0}^{\;\!t}
\|\, (\uu,\bb)(\cdot,\tau) \,
\|_{\mbox{}_{\scriptstyle L^{\infty}(\mathbb{R}^{4})}}
\|\,  \,\!(D\uu,D\bb)(\cdot,\tau) \,
\|_{\mbox{}_{\scriptstyle L^{2}(\mathbb{R}^{4})}}
\|\,  \,\!(D^{2}\uu,D^{2}\bb)(\cdot,\tau) \,
\|_{\mbox{}_{\scriptstyle L^{2}(\mathbb{R}^{4})}}
d\tau
\\
\leq\;
\|\,  \,\!(D\uu,D\bb)(\cdot,t_0) \,
\|_{\mbox{}_{\scriptstyle L^{2}(\mathbb{R}^{4})}}
^{\:\!2}
\\ +\;\:\!
C
\!\!\,\!
\int_{\scriptstyle t_0}^{\;\!t}
\|\,  \,\!(D\uu,D\bb)(\cdot,\tau) \,
\|_{\mbox{}_{\scriptstyle L^{2}(\mathbb{R}^{4})}}
\|\,  \,\!(D^{2}\uu,D^{2}\bb)(\cdot,\tau) \,
\|_{\mbox{}_{\scriptstyle L^{2}(\mathbb{R}^{4})}}
d\tau,
\end{multline*}
where we have used the Sobolev-Nirenberg-Gagliardo (SNG) inequalties (see (2.12)). Now, proceeding as in the 3D case we get
\begin{equation*}
\|(D\uu,D\bb)(\cdot,t)\|_{L^2(\RR^4)} \leq \|(D\uu,D\bb)(\cdot,t_0)\|_{L^2(\RR^4)}.
\end{equation*}
and consequently as in (2.4) one has
\begin{equation}
\tag{2.5}
\lim_{t\,\rightarrow\,\infty}
\,
t^{\:\!1/2}
\,
\|\, (D\uu,D\bb)(\cdot,t) \,
\|_{\mbox{}_{\scriptstyle L^{2}(\mathbb{R}^{n})}}
=\; 0,
\quad \;\;\,
2 \leq n \leq 4.
\end{equation}
\par In order to derive some Sobolev inequalities, we observe, by (1.7$e$), that
\begin{equation*}
\tag{2.6}
\| \uu \|_{L^q(\RR^n)} \leq \| (\uu,\bb )\|_{L^q(\RR^n)},
\end{equation*}
for $ 1 \leq q \leq \infty$.
%By the equivalence of norms, one has
%\begin{equation*}
%\tag{2.1$b$}
%\| \uu \|_{L^p(\RR^n)} \leq C \| (\uu,\bb )\|_{L^p(\RR^n)},
%\end{equation*}
%for some constant $C > 0 $.
The study of Leray solutions in dimension $ n \leq 4 $
is facilitated by the fact that
they are necessarily smooth for large $t$.
A further simplification for
\mbox{$ \:\!n = 2, 3 \:\!$}
is that
pointwise values of functions can be
estimated in terms of $H^{2}\!$ norms
and so we begin
with this case.
One has \\
\mbox{} \vspace{-0.650cm} \\
\begin{equation}
\tag{2.7$a$}
\|\, u \,
\|_{\mbox{}_{\scriptstyle L^{\infty}(\mathbb{R}^{2})}}
\:\!\leq\;
\|\, u \,
\|_{\mbox{}_{\scriptstyle L^{2}(\mathbb{R}^{2})}}
^{\:\!1/2}
\:\!
\|\, D^{2} u \,
\|_{\mbox{}_{\scriptstyle L^{2}(\mathbb{R}^{2})}}
^{\:\!1/2}
\end{equation}
\mbox{} \vspace{-0.175cm} \\
for arbitrary $ \:\! u \in H^{2}(\mathbb{R}^{2}) $;
likewise, \\
\mbox{} \vspace{-0.650cm} \\
\begin{equation}
\tag{2.7$b$}
\|\, u \,
\|_{\mbox{}_{\scriptstyle L^{\infty}(\mathbb{R}^{3})}}
\:\!\leq\;
\|\, u \,
\|_{\mbox{}_{\scriptstyle L^{2}(\mathbb{R}^{3})}}
^{\:\!1/4}
\:\!
\|\, D^{2} u \,
\|_{\mbox{}_{\scriptstyle L^{2}(\mathbb{R}^{3})}}
^{\:\!3/4}
\end{equation}
for $ u \in H^{2}(\mathbb{R}^{3}) $.
\!These are easily shown
by Fourier transform and Parseval's identity \linebreak
(see e.g.\;\cite{SchutzZiebellZinganoZingano2015},
where the optimal versions of (2.7) and their
higher dimensional analogues are obtained.
By Fourier transform,
we also get
(for any $n$): \\
\mbox{} \vspace{-0.650cm} \\
\begin{equation}
\tag{2.8$a$}
\|\, D \:\! u  \,
\|_{\mbox{}_{\scriptstyle L^{2}(\mathbb{R}^{n})}}
\:\!\leq\;
\|\, u \,
\|_{\mbox{}_{\scriptstyle L^{2}(\mathbb{R}^{n})}}
^{\:\!1/2}
\,\!
\|\, D^{2} u \,
\|_{\mbox{}_{\scriptstyle L^{2}(\mathbb{R}^{n})}}
^{\:\!1/2}
\end{equation}
or,
more generally, \\
\mbox{} \vspace{-0.600cm} \\
\begin{equation}
\tag{2.8$b$}
\|\, D^{\ell} \,\!u \,
\|_{\mbox{}_{\scriptstyle L^{2}(\mathbb{R}^{n})}}
\:\!\leq\;
\|\, u \,
\|_{\mbox{}_{\scriptstyle L^{2}(\mathbb{R}^{n})}}
^{\:\!1 - \theta}
\,\!
\|\, D^{m} \,\! u \,
\|_{\mbox{}_{\scriptstyle L^{2}(\mathbb{R}^{n})}}
^{\:\!\theta}
\!\:\!,
\quad \;\;
\theta \,=\:
\mbox{\small $ {\displaystyle \frac{\ell}{m} }$}
\end{equation}
\!Combining (2.6), (2.7) and (2.8),
we get the following
basic inequalities. \\
\\
\textbf{Lemma 1.}
\textit{For $n=2$, one has}
\begin{equation*}
\tag{2.9$a$}
\|\, (\uu,\bb) \,
\|_{\mbox{}_{\scriptstyle L^{\infty}(\mathbb{R}^{2})}}
\,\!
\|\,  \:\!(D\uu,D\bb) \,
\|_{\mbox{}_{\scriptstyle L^{2}(\mathbb{R}^{2})}}
\:\!\leq\; C
\|\, (\uu,\bb) \,
\|_{\mbox{}_{\scriptstyle L^{2}(\mathbb{R}^{2})}}
\:\!
\|\,  \,\!(D^2\uu,D^2\bb) \,
\|_{\mbox{}_{\scriptstyle L^{2}(\mathbb{R}^{2})}}
\!\:\!,
\end{equation*}
\mbox{} \vspace{-0.800cm} \\
\begin{equation*}
\tag{2.9$b$}
\|\, (\uu,\bb) \,
\|_{\mbox{}_{\scriptstyle L^{\infty}(\mathbb{R}^{2})}}
\,\!
\|\,  \,\!(D^2\uu,D^2\bb) \,
\|_{\mbox{}_{\scriptstyle L^{2}(\mathbb{R}^{2})}}
\:\!\leq\; C
\|\, (\uu,\bb) \,
\|_{\mbox{}_{\scriptstyle L^{2}(\mathbb{R}^{2})}}
\:\!
\|\,  \,\!(D^3\uu,D^3\bb) \,
\|_{\mbox{}_{\scriptstyle L^{2}(\mathbb{R}^{2})}}
\!\:\!,
\end{equation*}
\mbox{} \vspace{-0.800cm} \\
\begin{equation*}
\tag{2.9$c$}
\begin{split}
\|\,  \:\!(D\uu,D\bb) \,
\|_{\mbox{}_{\scriptstyle L^{\infty}(\mathbb{R}^{2})}}
\,\!
\|\,  \:\!(D\uu,D\bb) \,
\|_{\mbox{}_{\scriptstyle L^{2}(\mathbb{R}^{2})}}
\:\! \\ \leq\; C
\|\, (\uu,\bb) \,
\|_{\mbox{}_{\scriptstyle L^{2}(\mathbb{R}^{2})}}
\:\!
\|\,  \,\!(D^3\uu,D^3\bb) \,
\|_{\mbox{}_{\scriptstyle L^{2}(\mathbb{R}^{2})}}
\!\:\!,
\end{split}
\end{equation*}
\mbox{} \vspace{-0.150cm} \\
\textit{and,
for general}
$ \:\!m \geq 2 $, $ 0 \leq \ell \leq m - 2 $\/: \\
\mbox{} \vspace{-0.625cm} \\
\begin{equation*}
\tag{2.9$d$}
\begin{split}
\|\,  \,\!(D^{\ell}\uu,D^{\ell}\bb) \,
\|_{\mbox{}_{\scriptstyle L^{\infty}(\mathbb{R}^{2})}}
\,\!
\|\,  \,\!(D^{m - \ell}\uu,D^{m - \ell}\bb) \,
\|_{\mbox{}_{\scriptstyle L^{2}(\mathbb{R}^{2})}}
\:\!\\ \leq\; C
\|\, (\uu,\bb) \,
\|_{\mbox{}_{\scriptstyle L^{2}(\mathbb{R}^{2})}}
\:\!
\|\,  \,\!(D^{m + 1}\uu,D^{m + 1}\bb) \,
\|_{\mbox{}_{\scriptstyle L^{2}(\mathbb{R}^{2})}}
\!\:\!,
\end{split}
\end{equation*}
\textit{for some $C>0$}.
\\
\\
Similarly,
in dimension $ \:\!n = 3 $.
\\
\textbf{Lemma 2.}
\textit{For $n=3$, one has}
\begin{equation*}
\tag{2.10$a$}
\begin{split}
\|\, (\uu,\bb) \,
\|_{\mbox{}_{\scriptstyle L^{\infty}(\mathbb{R}^{3})}}
\,\!
\|\,  \:\!(D\uu,D\bb) \,
\|_{\mbox{}_{\scriptstyle L^{2}(\mathbb{R}^{3})}}
\:\!\\ \leq\; C
\|\, (\uu,\bb) \,
\|_{\mbox{}_{\scriptstyle L^{2}(\mathbb{R}^{3})}}
^{\:\!1/2}
\:\!
\|\,  \:\!(D\uu,D\bb) \,
\|_{\mbox{}_{\scriptstyle L^{2}(\mathbb{R}^{3})}}
^{\:\!1/2}
\:\!
\|\,  \,\!(D^{2}\uu,D^{2}\bb) \,
\|_{\mbox{}_{\scriptstyle L^{2}(\mathbb{R}^{3})}}
\!\:\!,
\end{split}
\end{equation*}
% \\
%
\begin{equation*}
\tag{2.10$b$}
\begin{split}
\|\, (\uu,\bb) \,
\|_{\mbox{}_{\scriptstyle L^{\infty}(\mathbb{R}^{3})}}
\,\!
\|\,  \,\!(D^{2}\uu,D^{2}\bb) \,
\|_{\mbox{}_{\scriptstyle L^{2}(\mathbb{R}^{3})}}
\:\! \\ \leq\;C
\|\, (\uu,\bb) \,
\|_{\mbox{}_{\scriptstyle L^{2}(\mathbb{R}^{3})}}
^{\:\!1/2}
\:\!
\|\,  \:\!(D\uu,D\bb) \,
\|_{\mbox{}_{\scriptstyle L^{2}(\mathbb{R}^{3})}}
^{\:\!1/2}
\:\!
\|\,  \,\!(D^{3}\uu,D^{3}\bb) \,
\|_{\mbox{}_{\scriptstyle L^{2}(\mathbb{R}^{3})}}
\!\:\!,
\end{split}
\end{equation*}
\\
\begin{equation*}
\tag{2.10$c$}
\begin{split}
\|\,  \:\!(D\uu,D\bb) \,
\|_{\mbox{}_{\scriptstyle L^{\infty}(\mathbb{R}^{3})}}
\,\!
\|\,  \:\!(D\uu,D\bb) \,
\|_{\mbox{}_{\scriptstyle L^{2}(\mathbb{R}^{3})}}
\:\! \\ \leq\;C
\|\, (\uu,\bb) \,
\|_{\mbox{}_{\scriptstyle L^{2}(\mathbb{R}^{3})}}
^{\:\!1/2}
\:\!
\|\,  \:\!(D\uu,D\bb) \,
\|_{\mbox{}_{\scriptstyle L^{2}(\mathbb{R}^{3})}}
^{\:\!1/2}
\:\!
\|\,  \,\!(D^{3}\uu,D^{3}\bb) \,
\|_{\mbox{}_{\scriptstyle L^{2}(\mathbb{R}^{3})}}
\!\:\!,
\end{split}
\end{equation*}
\\
\begin{equation*}
\tag{2.10$d$}
\begin{split}
\|\,  \:\!(D\uu,D\bb) \,
\|_{\mbox{}_{\scriptstyle L^{\infty}(\mathbb{R}^{3})}}
\,\!
\|\,  \,\!(D^{2}\uu,D^{2}\bb) \,
\|_{\mbox{}_{\scriptstyle L^{2}(\mathbb{R}^{3})}}
\:\!\\ \leq\;C
\|\, (\uu,\bb) \,
\|_{\mbox{}_{\scriptstyle L^{2}(\mathbb{R}^{3})}}
^{\:\!3/4}
\:\!
\|\,  \,\!(D^{2}\uu,D^{2}\bb) \,
\|_{\mbox{}_{\scriptstyle L^{2}(\mathbb{R}^{3})}}
^{\:\!1/4}
\:\!
\|\,  \,\!(D^{4}\uu,D^{4}\bb) \,
\|_{\mbox{}_{\scriptstyle L^{2}(\mathbb{R}^{3})}}
\!\:\!,
\end{split}
\end{equation*}
 \\
\begin{equation*}
\tag{2.10$e$}
\begin{split}
\|\,  \,\!(D^{2}\uu,D^{2}\bb) \,
\|_{\mbox{}_{\scriptstyle L^{\infty}(\mathbb{R}^{3})}}
\,\!
\|\,  \,\!(D^{2}\uu,D^{2}\bb) \,
\|_{\mbox{}_{\scriptstyle L^{2}(\mathbb{R}^{3})}}
\:\!\\ \leq\;C
\|\, (\uu,\bb) \,
\|_{\mbox{}_{\scriptstyle L^{2}(\mathbb{R}^{3})}}
^{\:\!3/4}
\:\!
\|\,  \,\!(D^{2}\uu,D^{2}\bb) \,
\|_{\mbox{}_{\scriptstyle L^{2}(\mathbb{R}^{3})}}
^{\:\!1/4}
\:\!
\|\,  \,\!(D^{5}\uu,D^{5}\bb) \,
\|_{\mbox{}_{\scriptstyle L^{2}(\mathbb{R}^{3})}}
\!\:\!,
\end{split}
\end{equation*}
\mbox{} \vspace{-0.150cm} \\
\\
\textit{and,
for general}
$ \:\!m \geq 3 $, $ 0 \leq \ell \leq m - 3 $\/: \\

\mbox{} \vspace{-0.700cm} \\
\begin{equation*}
\tag{2.10$f$}
\begin{split}
\|\,  \,\!(D^{\ell}\uu,D^{\ell}\bb) \,
\|_{\mbox{}_{\scriptstyle L^{\infty}(\mathbb{R}^{3})}}
\,\!
\|\,  \,\!(D^{m - \ell}\uu,D^{m - \ell}\bb) \,
\|_{\mbox{}_{\scriptstyle L^{2}(\mathbb{R}^{3})}}
\:\! \\ \leq\;C
\|\, (\uu,\bb) \,
\|_{\mbox{}_{\mbox{}_{\scriptstyle L^{2}(\mathbb{R}^{3})}}}
^{{\scriptstyle \:\!
		\frac{\scriptstyle \ell \,+\, 3/2}{\scriptstyle \ell \,+\, 2} }}
\,
\|\,  \,\!(D^{\ell + 2}\uu,D^{\ell + 2}\bb) \,
\|_{\mbox{}_{\mbox{}_{\scriptstyle L^{2}(\mathbb{R}^{3})}}}
^{{\scriptstyle \:\!
		\frac{\scriptstyle 1/2}{\scriptstyle \ell \,+\, 2} }}
\|\,  \,\!(D^{m + 1}\uu,D^{m + 1}\bb) \,
\|_{\mbox{}_{\scriptstyle L^{2}(\mathbb{R}^{3})}}
\!\;\!,
\end{split}
\end{equation*}
\textit{for some $ C>0$}.
\\
In dimension $n=4$, we start with the fundamental
Sobolev inequality,
\begin{equation*}
\tag{2.11}
\| u \|_{L^4(\RR^4)} \leq \| Du\|_{L^2(\RR^4)}.
\end{equation*}
Hence, using (2.6), (2.8) and (2.11), we have the result below.
\\
\\
\textbf{Lemma 3.}
\textit{For all $ \:\! m \geq 1 $,
	$ 0 \leq \ell \leq m - 1 $, one actually has}
\begin{equation*}
\tag{2.12}
\begin{split}
\|\,  \,\!(D^{\ell}\uu,D^{\ell}\bb)\,
\|_{\mbox{}_{\scriptstyle L^{4}(\mathbb{R}^{4})}}
\|\,  \,\!(D^{m - \ell}\uu,D^{m - \ell}\bb)\,
\|_{\mbox{}_{\scriptstyle L^{4}(\mathbb{R}^{4})}}
\!\;\!
\\
\leq\; C
\|\,  \:\!(D\uu,D\bb)\,
\|_{\mbox{}_{\scriptstyle L^{2}(\mathbb{R}^{4})}}
\|\,  \,\!(D^{m + 1}\uu,D^{m + 1}\bb)\,
\|_{\mbox{}_{\scriptstyle L^{2}(\mathbb{R}^{4})}},
\end{split}
\end{equation*}
\textit{for some $C>0$}.
\par When we derive energy inequalities for higher order derivatives
of Leray solutions to MHD equations, the importance of lemmas above becomes clear. In euclidean plane $\RR^2$, it turns out that all solutions of MHD system (1.1)  are the same, i.e., the uniqueness is well established; the  solutions are also to be smooth, in other words,
$ (\uu,\bb) \in C^\infty(\RR^2 \times (0,\infty)) $ and moreover
$ {\displaystyle
	\;\!
(	\uu,\bb)(\cdot,t)
	\in
	C((\:\!0, \infty), \:\!\mbox{\boldmath $H$}^{\!\;\!m}\!\;\!(\mathbb{R}^{2}))
} $
for all $ m \geq 0 $.
When $ n > 2 $,
the absence of smoothness previously cited
complicates the study
of Leray solutions;
in particular,
their uniqueness and precise regularity
properties
are still unresolved as in the Navier-Stokes system case.
\par Now, we will generalize the argument in
\cite{SchutzZinganoZingano2015}
for the MHD system (1.1) in dimension $n=4$.
Since $ \mbox{\boldmath $u$}(\cdot,t) $
is smooth
for large $ \:\!t $,
it can be written as \\
\mbox{} \vspace{-0.700cm} \\
\begin{equation}
\tag{2.13}
\mbox{} \hspace{+0.500cm}
\mbox{\boldmath $u$}(\cdot,t)
\;=\;
e^{\:\!\mu\:\!\Delta (t - t_0)}
\:\!
\mbox{\boldmath $u$}(\cdot,t_0)
\,-
\int_{\!\:\!t_0}^{\;\!t}
\!\!\:\!
e^{\:\!\mu\:\!\Delta (t - \tau)}
\:\!
\mbox{\boldmath $Q_1$}(\cdot,\tau)
\,
d\tau,
\quad \;\;\,
t > t_0
\end{equation}
\mbox{} \vspace{-0.225cm} \\
for $ \:\!t_0 \!\;\!$
large enough,
where
$ {\displaystyle
	\:\!
	\mbox{\boldmath $Q_1$} \!\;\!=\,\!
	\mbox{\boldmath $u$} \!\;\!\cdot\!\;\!\nabla \mbox{\boldmath $u$}
	\:\!+ \nabla P - \bb \cdot \nabla \bb
} $,
and
$ \:\!e^{\:\!\mu\:\!\Delta t} \!\!\;\!\;\!$
denotes the heat semigroup.
%
%
% -------------------------- Proof for n = 4:
%
%
From (2.13),
we get
\begin{equation*}
\begin{split}
	\|\, \mbox{\boldmath $u$}(\cdot,t) \,
	\|_{\mbox{}_{\scriptstyle L^{2}(\mathbb{R}^{4})}}
	\leq\;
	\|\, \mbox{\boldmath $v$}_{0}(\cdot,t) \,
	\|_{\mbox{}_{\scriptstyle L^{2}(\mathbb{R}^{4})}}
	\:\!+
	\int_{\!\:\!t_0}^{\;\!t}
	\!\!\:\!
	\|\, \mbox{\boldmath $u$}(\cdot,\tau) \!\;\!\cdot\!\;\!
	\nabla \mbox{\boldmath $u$}(\cdot,\tau) \,
	\|_{\mbox{}_{\scriptstyle L^{2}(\mathbb{R}^{4})}}
	\;\!
	d\tau
	\\ +
	\int_{\!\:\!t_0}^{\;\!t}
	\!\!\:\!
	\|\, \mbox{\boldmath $b$}(\cdot,\tau) \!\;\!\cdot\!\;\!
	\nabla \mbox{\boldmath $b$}(\cdot,\tau) \,
	\|_{\mbox{}_{\scriptstyle L^{2}(\mathbb{R}^{4})}}
	\;\!
	d\tau
\\
	\leq\;
	\|\, \mbox{\boldmath $v$}_{0}(\cdot,t) \,
	\|_{\mbox{}_{\scriptstyle L^{2}(\mathbb{R}^{4})}}
	\:\!+\,
	\sqrt{\:\!2\;\!\:\!} \!\!
	\int_{\!\:\!t_0}^{\;\!t}
	\!\!\:\!
	\|\, \mbox{\boldmath $u$}(\cdot,\tau) \,
	\|_{\mbox{}_{\scriptstyle L^{4}(\mathbb{R}^{4})}}
	\,\!
	\|\, D \:\! \mbox{\boldmath $u$}(\cdot,\tau) \,
	\|_{\mbox{}_{\scriptstyle L^{4}(\mathbb{R}^{4})}}
	\;\!
	d\tau
	\\
	\sqrt{\:\!2\;\!\:\!} \!\!
	\int_{\!\:\!t_0}^{\;\!t}
	\!\!\:\!
	\|\, \mbox{\boldmath $b$}(\cdot,\tau) \,
	\|_{\mbox{}_{\scriptstyle L^{4}(\mathbb{R}^{4})}}
	\,\!
	\|\, D \:\! \mbox{\boldmath $b$}(\cdot,\tau) \,
	\|_{\mbox{}_{\scriptstyle L^{4}(\mathbb{R}^{4})}}
	\;\!
	d\tau \\
	\leq\;
	\|\, \mbox{\boldmath $v$}_{0}(\cdot,t) \,
	\|_{\mbox{}_{\scriptstyle L^{2}(\mathbb{R}^{4})}}
	\:\!+\,
	2\sqrt{\:\!2\;\!\:\!} \!\!
	\int_{\!\:\!t_0}^{\;\!t}
	\!\!\:\!
	\|\, (D \:\!\mbox{\boldmath $u$},D \:\!\mbox{\boldmath $b$})(\cdot,\tau) \,
	\|_{\mbox{}_{\scriptstyle L^{2}(\mathbb{R}^{4})}}
	\,\!
	\|\, (D^{2} \,\! \mbox{\boldmath $u$},D^{2} \,\! \mbox{\boldmath $b$})(\cdot,\tau) \,
	\|_{\mbox{}_{\scriptstyle L^{2}(\mathbb{R}^{4})}}
	\;\!
	d\tau
\end{split}
\end{equation*}
by (2.11),
where
$ {\displaystyle
	\:\!
	\mbox{\boldmath $v$}_{0}(\cdot,t)
	\!:=
	e^{\:\!\mu\:\!\Delta (t - t_0)}
	\:\!
	\mbox{\boldmath $u$}(\cdot,t_0)
} $,
and using that
(by Helmholtz projection
or directly by Fourier transform
\cite{KreissHagstromLorenzZingano2003}):
$ {\displaystyle
	\!\;\!
	\|\, \mbox{\boldmath $Q_1$}(\cdot,\tau) \,
	\|_{\scriptstyle L^{2}(\mathbb{R}^{n})}
	\!\leq
	\|\, \mbox{\boldmath $u$}(\cdot,\tau) \!\:\!\cdot\!\:\!
	\nabla \mbox{\boldmath $u$}(\cdot,\tau) \,
	\|_{\scriptstyle L^{2}(\mathbb{R}^{n})}
	\!\:\!
} $.
% (for any $n$).
This shows that, given $ \:\!\epsilon > 0 $,
taking $ \:\!t_0 \gg 1 $
we get
$ {\displaystyle
	\|\, \mbox{\boldmath $u$}(\cdot,t) \,
	\|_{\scriptstyle L^{2}(\mathbb{R}^{4})}
	\!\;\!< \epsilon
	\;\!
} $
for all $ t $ large enough,
since
the integrand on the
righthand side above is in
$ L^{1}((\:\!t_{\ast}, \infty)) $.
One can repeat the previous analysis for $\bb(\cdot,t)$ using
$ {\displaystyle
	\:\!
	\mbox{\boldmath $Q_2$} \!\;\!=\,\!
	\mbox{\boldmath $u$} \!\;\!\cdot\!\;\!\nabla \mbox{\boldmath $b$}
	 - \bb \cdot \nabla \uu
} $
 and obtain
\begin{equation*}
\| (\uu,\bb)(\cdot,t) \|_{L^2(\RR^4)} \to 0, \text{     as   } t \to \infty,
\end{equation*}
which implies (with (1.4)) that
\begin{equation*}
\tag{2.14}
\| (\uu,\bb)(\cdot,t) \|_{L^2(\RR^n)} \to 0, \text{     as   } t \to \infty,
\end{equation*}
for $ n=2,3,4$.
%
%
% ***************************************************
% *                                                 *
% *                  Section 3                      *
% *                                                 *
% *                                                 *
% ***************************************************
%
\\
\par
{\bf 3. Proof of Main Theorem }
\\
\par Let $(\uu,\bb)(\cdot,t)$ be any given Leray solution of the system (1.1). Observe that, by (2.14) and (2.5), the result is true for $s = 0$ and $s = 1$. Our strategy will be to show that the main theorem is valid for $s > 0$ integer, i.e.,
\begin{equation*}
\tag{3.1}
\begin{split}
\lim_{t \to \infty} t^{m/2} \| (D^m\uu,D^m\bb)(\cdot,t) \|_{{L^2(\RR^n)}} = 0,
\end{split}
\end{equation*}
for all $m \geq 0$ integer.
\par By (2.5) and (2.14), given $\epsilon >0$, there exist $t_0 > t_* $ (see (1.3))
sufficiently large such as
\begin{equation*}
\tag{(3.2$a$)}
\|(\uu,\bb)(\cdot,t)\|_{L^2(\RR^n)} \leq \epsilon
\end{equation*}
and
\begin{equation*}
\tag{3.2$b$}
t^{1/2}\|(D\uu,D\bb)(\cdot,t)\|_{L^2(\RR^n)} \leq \epsilon,
\end{equation*}

for all $t \geq t_0$.
%-----------------------------------------------------
%
%                         m = 1
%
%------------------------------------------------------

\par Starting with $n =3$, let $t_* \geq 0$ be chosen so that (1.3) holds. Now, Differentiating (1.1$a$) and (1.1$b$) with respect to $ x_\ell$, taking the dot product of (1.1$a$) and (1.1$b$) by $(t-t_0) D_\ell \uu$ and $(t-t_0) D_\ell \bb$, respectively, and integrating the result on $\RR^3 \times [t_0,t] $, we get summing over $ 1 \leq \ell \leq 3$,

\begin{multline*}
(\:\!t - t_{0})
\,
\|\,  \,\!(D\uu,D\bb)(\cdot,t) \,
\|_{\mbox{}_{\scriptstyle L^{2}(\mathbb{R}^{3})}}
^{\:\!2}
\!\;\! \\ +\:
2\, \min\{\mu,\nu\} \!\!\;\!
\int_{\scriptstyle t_0}^{\;\!t}
\!\!\:\!
(\tau - t_{0})
\,
\|\,  \,\!(D^{2}\uu,D^{2}\bb) (\cdot,\tau) \,
\|_{\mbox{}_{\scriptstyle L^{2}(\mathbb{R}^{3})}}
^{\:\!2}
d\tau
\\
\leq\;
\!\!\;\!
\int_{\scriptstyle t_0}^{\;\!t}
\!\!\;\!
\,
\|\,  \,\!(D\uu,D\bb)(\cdot,\tau) \,
\|_{\mbox{}_{\scriptstyle L^{2}(\mathbb{R}^{3})}}
^{\:\!2}
d\tau
\\ +\;\:\!
C
\!\!\,\!
\int_{\scriptstyle t_0}^{\;\!t}
\!\!\;\!
(\tau - t_{0})
\|\, (\uu,\bb)(\cdot,\tau) \,
\|_{\mbox{}_{\scriptstyle L^{\infty}(\mathbb{R}^{3})}}
\|\,  \,\!(D\uu,D\bb)(\cdot,\tau) \,
\|_{\mbox{}_{\scriptstyle L^{2}(\mathbb{R}^{3})}}
\|\,  \,\!(D^{2}\uu,D^{2}\bb)(\cdot,\tau) \,
\|_{\mbox{}_{\scriptstyle L^{2}(\mathbb{R}^{3})}}
d\tau
\\
\leq\;
\!\!\;\!
\int_{\scriptstyle t_0}^{\;\!t}
\!\!\;\!
\,
\|\,  \,\!(D\uu,D\bb)(\cdot,\tau) \,
\|_{\mbox{}_{\scriptstyle L^{2}(\mathbb{R}^{3})}}
^{\:\!2}
d\tau
\\ +\;\:\!
C
\!\!\,\!
\int_{\scriptstyle t_0}^{\;\!t}
\!\!\;\!
(\tau - t_{0})
\|\, (\uu,\bb)(\cdot,\tau) \,
\|^{1/2}_{\mbox{}_{\scriptstyle L^{2}(\mathbb{R}^{3})}}
\|\,  \,\!(D\uu,D\bb)(\cdot,\tau) \,
\|^{1/2}_{\mbox{}_{\scriptstyle L^{2}(\mathbb{R}^{3})}}
\|\,  \,\!(D^{2}\uu,D^{2}\bb)(\cdot,\tau) \,
\|_{\mbox{}_{\scriptstyle L^{2}(\mathbb{R}^{3})}}^{2}
d\tau,
\end{multline*}
where we have used integration by parts, (1.1$c$) and (2.10$a$). Therefore  by (2.5) and (2.14), for $t _0>t_*$ sufficiently large, we have,
\begin{equation*}
\tag{3.3}
\begin{split}
(\:\!t - t_{0})
\,
\|\,  \,\!(D\uu,D\bb)(\cdot,t) \,
\|_{\mbox{}_{\scriptstyle L^{2}(\mathbb{R}^{3})}}
^{\:\!2}
\!\;\!  +\:
C \!\!\;\!
\int_{\scriptstyle t_0}^{\;\!t}
\!\!\:\!
(\tau - t_{0})
\,
\|\,  \,\!(D^{2}\uu,D^{2}\bb) (\cdot,\tau) \,
\|_{\mbox{}_{\scriptstyle L^{2}(\mathbb{R}^{3})}}
^{\:\!2}
d\tau
\\
\leq\;
\!\!\;\!
\int_{\scriptstyle t_0}^{\;\!t}
\!\!\;\!
\,
\|\,  \,\!(D\uu,D\bb)(\cdot,\tau) \,
\|_{\mbox{}_{\scriptstyle L^{2}(\mathbb{R}^{3})}}
^{\:\!2}
d\tau.
\end{split}
\end{equation*}
For some constant $C>0$.

%---------------------------------------------------------
%
%
%                       m = 2
%
%
%---------------------------------------------------------

\par Now, we go to the next step similarly: differentiating (1.1$a$) and
(1.1$b$) twice (with respect to $x_{\ell_1}$, $x_{\ell_2}$, for example), multiplying  (1.1$a$) and (1.1$b$) by $(t - t_0)^2 D_{\mbox{}_{\scriptstyle \!\:\!\ell_{1}}}
\! D_{\mbox{}_{\scriptstyle \!\:\!\ell_{2}}}
\!\mbox{\boldmath $u$}(x,t) $ and by $ D_{\mbox{}_{\scriptstyle \!\:\!\ell_{1}}}
\! D_{\mbox{}_{\scriptstyle \!\:\!\ell_{2}}}
\!\mbox{\boldmath $b$}(x,t) $, respectively,
we get, integrating the result
on $ \mathbb{R}^{3} \!\:\!\times\!\;\![\,t_0, \:\!t\;\!] $,
$ \:\! t \geq t_0 $ and summing over $1\leq\ell_1,\ell_2 \leq 3$,
\begin{multline*}
(\:\!t - t_{0})^{2}
\,
\|\,  \,\!(D^{2}\uu,D^{2}\bb)(\cdot,t) \,
\|_{\mbox{}_{\scriptstyle L^{2}(\mathbb{R}^{3})}}
^{\:\!2}
\!\;\! \\ +\:
2\, \min\{\mu,\nu\} \!\!\;\!
\int_{\scriptstyle t_0}^{\;\!t}
\!\!\:\!
(\tau - t_{0})^{2}
\,
\|\,  \,\!(D^{3}\uu,D^{3}\bb) (\cdot,\tau) \,
\|_{\mbox{}_{\scriptstyle L^{2}(\mathbb{R}^{3})}}
^{\:\!2}
d\tau
\\
\leq\;
2
\!\!\;\!
\int_{\scriptstyle t_0}^{\;\!t}
\!\!\;\!
(\tau - t_{0})
\,
\|\, \,\!( D^{2}\uu, D^{2}\bb)(\cdot,\tau) \,
\|_{\mbox{}_{\scriptstyle L^{2}(\mathbb{R}^{3})}}
^{\:\!2}
d\tau
\\ +\;\:\!
C
\!\!\,\!
\int_{\scriptstyle t_0}^{\;\!t}
\!\!\;\!
(\tau - t_{0})^{2}
\,
\Bigl\{\,
\|\, (\uu,\bb)(\cdot,\tau) \,
\|_{\mbox{}_{\scriptstyle L^{\infty}(\mathbb{R}^{3})}}
\|\,  \,\!(D^{2}\uu,D^{2}\bb)(\cdot,\tau) \,
\|_{\mbox{}_{\scriptstyle L^{2}(\mathbb{R}^{3})}}
\\
+\;\;\!
\|\, \:\!(D\uu,D\bb)(\cdot,\tau) \,
\|_{\mbox{}_{\scriptstyle L^{\infty}(\mathbb{R}^{3})}}
\|\,  (D\uu,D\bb)(\cdot,\tau) \,
\|_{\mbox{}_{\scriptstyle L^{2}(\mathbb{R}^{3})}}
\;\!\Bigr\}
\,
\|\,  \,\!(D^{3}\uu,D^{3}\bb)(\cdot,\tau) \,
\|_{\mbox{}_{\scriptstyle L^{2}(\mathbb{R}^{3})}}
d\tau.
\end{multline*}
Using (2.10$b$) and (2.10$c$), we have,
\begin{multline*}
(\:\!t - t_{0})^{2}
\,
\|\,  \,\!(D^{2}\uu,D^{2}\bb)(\cdot,t) \,
\|_{\mbox{}_{\scriptstyle L^{2}(\mathbb{R}^{3})}}
^{\:\!2}
\!\;\! \\ +\:
2\, \min\{\mu,\nu\} \!\!\;\!
\int_{\scriptstyle t_0}^{\;\!t}
\!\!\:\!
(\tau - t_{0})^{2}
\,
\|\,  \,\!(D^{3}\uu,D^{3}\bb) (\cdot,\tau) \,
\|_{\mbox{}_{\scriptstyle L^{2}(\mathbb{R}^{3})}}
^{\:\!2}
d\tau
\\
\leq\;
2
\!\!\;\!
\int_{\scriptstyle t_0}^{\;\!t}
\!\!\;\!
(\tau - t_{0})
\,
\|\, D^{2} \,\!(\uu,\bb)(\cdot,\tau) \,
\|_{\mbox{}_{\scriptstyle L^{2}(\mathbb{R}^{3})}}
^{\:\!2}
d\tau
\\ +\;\:\!
C
\!\!\,\!
\int_{\scriptstyle t_0}^{\;\!t}
\!\!\;\!
(\tau - t_{0})^{2}
\|\, (\uu,\bb)(\cdot,\tau) \,
\|^{1/2}_{\mbox{}_{\scriptstyle L^{2}(\mathbb{R}^{3})}}
\|\,  \,\!(D\uu,D\bb)(\cdot,\tau) \,
\|^{1/2}_{\mbox{}_{\scriptstyle L^{2}(\mathbb{R}^{3})}}
\|\,  \,\!(D^{3}\uu,D^{3}\bb)(\cdot,\tau) \,
\|_{\mbox{}_{\scriptstyle L^{2}(\mathbb{R}^{3})}}^{2}
d\tau,
\end{multline*}
for some constant $ C>0 $
(whose value need not concern us). Hence, by (2.5) and (2.14), there exist
$t_0 > t_*$ sufficiently large such as,
\begin{equation*}
\tag{3.4}
\begin{split}
(\:\!t - t_{0})^2
\,
\|\,  \,\!(D^2\uu,D^2\bb)(\cdot,t) \,
\|_{\mbox{}_{\scriptstyle L^{2}(\mathbb{R}^{3})}}
^{\:\!2}
\!\;\!  +\:
C \!\!\;\!
\int_{\scriptstyle t_0}^{\;\!t}
\!\!\:\!
(\tau - t_{0})^2
\,
\|\,  \,\!(D^{3}\uu,D^{3}\bb) (\cdot,\tau) \,
\|_{\mbox{}_{\scriptstyle L^{2}(\mathbb{R}^{3})}}
^{\:\!2}
d\tau
\\
\leq\;
2
\!\!\;\!
\int_{\scriptstyle t_0}^{\;\!t}
\!\!\;\!
(\tau - t_0)
\|\,  \,\!(D^2\uu,D^2\bb)(\cdot,\tau) \,
\|_{\mbox{}_{\scriptstyle L^{2}(\mathbb{R}^{3})}}
^{\:\!2}
d\tau.
\end{split}
\end{equation*}
%----------------------------------------------------------------
%
%
%                     proof of m = 2 property
%
%----------------------------------------------------------------
Given $\epsilon > 0$ arbitrary, there exist $t_0 > t_*$ large enough so that,
\begin{equation*}
\begin{split}
\int_{\scriptstyle t_0}^{\;\!t}
\!\!\;\!
\,
\|\,  \,\!(D\uu,D\bb)(\cdot,\tau) \,
\|_{\mbox{}_{\scriptstyle L^{2}(\mathbb{R}^{3})}}
^{\:\!2}
d\tau \leq \epsilon,
\end{split}
 \end{equation*}
by the fundamental energy inequality (1.2) and (3.2). Hence, by (3.3), one has,
\begin{equation*}
\begin{split}
\int_{\scriptstyle t_0}^{\;\!t}
\!\!\:\!
(\tau - t_{0})
\,
\|\,  \,\!(D^{2}\uu,D^{2}\bb) (\cdot,\tau) \,
\|_{\mbox{}_{\scriptstyle L^{2}(\mathbb{R}^{3})}}
^{\:\!2}
d\tau \leq \epsilon.
\end{split}
\end{equation*}
Using (3.4), we conclude that,
\begin{equation*}
\begin{split}
t^2 \Bigg(\frac{\:\!t - t_{0}}{t}\Bigg)^2
\,
\|\,  \,\!(D^2\uu,D^2\bb)(\cdot,t) \,
\|_{\mbox{}_{\scriptstyle L^{2}(\mathbb{R}^{3})}}
^{\:\!2} = (\:\!t - t_{0})^2
\,
\|\,  \,\!(D^2\uu,D^2\bb)(\cdot,t) \,
\|_{\mbox{}_{\scriptstyle L^{2}(\mathbb{R}^{3})}}
^{\:\!2} \leq \epsilon.
\end{split}
\end{equation*}
Consequently
\begin{equation*}
\tag{3.5}
\begin{split}
t \|\,  \,\!(D^2\uu,D^2\bb)(\cdot,t) \,
\|_{\mbox{}_{\scriptstyle L^{2}(\mathbb{R}^{3})}} \to 0,
\text{   as  } t \to \infty
\end{split}
\end{equation*}
and (3.5) solves (3.1) for $m = 2$.
Similarly, we go to the next step and use the previous decay (3.5) and the Sobolev inequalities (2.10$e$) and (2.10$f$)  to obtain the 3rd order decay. Now, by induction, we have,
\begin{equation*}
\begin{split}
(\:\!t - t_{0})^{m}
\,
\|\,  \,\!(D^{m}\uu,D^{m}\bb)(\cdot,t) \,
\|_{\mbox{}_{\scriptstyle L^{2}(\mathbb{R}^{3})}}
^{\:\!2}
\!\;\! \\ +\:
2\, \min\{\mu,\nu\} \!\!\;\!
\int_{\scriptstyle t_0}^{\;\!t}
\!\!\:\!
(\tau - t_{0})^{m}
\,
\|\,  \,\!(D^{m+1}\uu,D^{m+1}\bb) (\cdot,\tau) \,
\|_{\mbox{}_{\scriptstyle L^{2}(\mathbb{R}^{3})}}
^{\:\!2}
d\tau
\\
\leq\;
m
\!\!\;\!
\int_{\scriptstyle t_0}^{\;\!t}
\!\!\;\!
(\tau - t_{0})^{m-1}
\,
\|\,  \,\!(D^{m}\uu,D^{m}\bb)(\cdot,\tau) \,
\|_{\mbox{}_{\scriptstyle L^{2}(\mathbb{R}^{3})}}
^{\:\!2}
d\tau
\\ +\;\:\!
C
\!\!\,\!
\int_{\scriptstyle t_0}^{\;\!t}
\!\!\;\!
(\tau - t_{0})^{m}
\|\,  \,\!(D^{m+1}\uu,D^{m+1}\bb)(\cdot,\tau) \,
\|_{\mbox{}_{\scriptstyle L^{2}(\mathbb{R}^{3})}}
\sum_{\ell\,=\,0}^{[\,m/2\,]}
\|\,  \,\! (D^{\ell}\uu,D^{\ell}\bb)(\cdot,\tau) \,
\|_{\mbox{}_{\scriptstyle L^{\infty}(\mathbb{R}^{3})}}
\\
\cdot \|\,  \,\!(D^{m - \ell}\uu,D^{m - \ell}\bb)(\cdot,\tau) \,
\|_{\mbox{}_{\scriptstyle L^{2}}(\mathbb{R}^{3})}
\;\!
d\tau,
\end{split}
\end{equation*}
for general $m \geq 3$, where $[m]$ is the integer part of $m$.
By (2.10$f$),
\begin{equation*}
\begin{split}
(\:\!t - t_{0})^m
\,
\|\,  \,\!(D^m\uu,D^m\bb)(\cdot,t) \,
\|_{\mbox{}_{\scriptstyle L^{2}(\mathbb{R}^{3})}}
^{\:\!2}
\!\;\!  +\:
C \!\!\;\!
\int_{\scriptstyle t_0}^{\;\!t}
\!\!\:\!
(\tau - t_{0})^m
\,
\|\,  \,\!(D^{m+1}\uu,D^{m+1}\bb) (\cdot,\tau) \,
\|_{\mbox{}_{\scriptstyle L^{2}(\mathbb{R}^{3})}}
^{\:\!2}
d\tau
\\
\leq\;
m
\!\!\;\!
\int_{\scriptstyle t_0}^{\;\!t}
\!\!\;\!
(\tau - t_0)^{m-1}
\|\,  \,\!(D^m\uu,D^m\bb)(\cdot,\tau) \,
\|_{\mbox{}_{\scriptstyle L^{2}(\mathbb{R}^{3})}}
^{\:\!2}
d\tau.
\end{split}
\end{equation*}
By the same previous argument, it follows that
\begin{equation*}
t^{m/2} \|(D^m\uu,D^m\bb)(\cdot,t)\|_{L^2(\RR^3)} \to 0 \text{  as }
t \to \infty,
\end{equation*}
which completes the proof of (3.1) for $n=3$. The proof of (3.1)
for $n=2$ is similar, using the inequalities of Lemma 1 (2.9) instead of Lemma 2.
%------------------------------------------------------------
%
%
%                       n = 4
%
%
%-------------------------------------------------------------

\par  We will now consider the $n = 4$ case. Basically, we will use the inequalities (2.11) and (2.12). However, the energy estimates will suffer some changes. So, let $(u,b)(\cdot,t) $ be any given Leray solution to (1.1). Differentiating (1.1$a$) and (1.1$b$) with respect to $ x_\ell$, taking the dot product of (1.1$a$) and (1.1$b$) by $(t-t_0) D_\ell \uu$ and $(t-t_0) D_\ell \bb$, respectively, and integrating the result on $\RR^4 \times [t_0,t] $, the energy estimate, summing over $ 1 \leq \ell \leq 4$, is now

\begin{multline*}
(\:\!t - t_{0})
\,
\|\,  \,\!(D\uu,D\bb)(\cdot,t) \,
\|_{\mbox{}_{\scriptstyle L^{2}(\mathbb{R}^{4})}}
^{\:\!2}
\!\;\! \\ +\:
2\, \min\{\mu,\nu\} \!\!\;\!
\int_{\scriptstyle t_0}^{\;\!t}
\!\!\:\!
(\tau - t_{0})
\,
\|\,  \,\!(D^{2}\uu,D^{2}\bb) (\cdot,\tau) \,
\|_{\mbox{}_{\scriptstyle L^{2}(\mathbb{R}^{4})}}
^{\:\!2}
d\tau
\\
\leq\;
\!\!\;\!
\int_{\scriptstyle t_0}^{\;\!t}
\!\!\;\!
\,
\|\,  \,\!(D\uu,D\bb)(\cdot,\tau) \,
\|_{\mbox{}_{\scriptstyle L^{2}(\mathbb{R}^{4})}}
^{\:\!2}
d\tau
\\ +\;\:\!
C
\!\!\,\!
\int_{\scriptstyle t_0}^{\;\!t}
\!\!\;\!
(\tau - t_{0})
\|\, (\uu,\bb)(\cdot,\tau) \,
\|_{\mbox{}_{\scriptstyle L^{4}(\mathbb{R}^{4})}}
\|\,  \,\!(D\uu,D\bb)(\cdot,\tau) \,
\|_{\mbox{}_{\scriptstyle L^{4}(\mathbb{R}^{4})}}
\|\,  \,\!(D^{2}\uu,D^{2}\bb)(\cdot,\tau) \,
\|_{\mbox{}_{\scriptstyle L^{2}(\mathbb{R}^{4})}}
d\tau
\\
\leq\;
\!\!\;\!
\int_{\scriptstyle t_0}^{\;\!t}
\!\!\;\!
\,
\|\,  \,\!(D\uu,D\bb)(\cdot,\tau) \,
\|_{\mbox{}_{\scriptstyle L^{2}(\mathbb{R}^{4})}}
^{\:\!2}
d\tau
\\ +\;\:\!
C
\!\!\,\!
\int_{\scriptstyle t_0}^{\;\!t}
\!\!\;\!
(\tau - t_{0})
\|\,  \,\!(D\uu,D\bb)(\cdot,\tau) \,
\|_{\mbox{}_{\scriptstyle L^{2}(\mathbb{R}^{4})}}
\|\,  \,\!(D^{2}\uu,D^{2}\bb)(\cdot,\tau) \,
\|^2_{\mbox{}_{\scriptstyle L^{2}(\mathbb{R}^{4})}}
d\tau,
\end{multline*}
by the H\"older inequality and (2.11) for vector field $(\uu,\bb)$.
Using (2.5), one has,
\begin{equation*}
\begin{split}
(\:\!t - t_{0})
\,
\|\,  \,\!(D\uu,D\bb)(\cdot,t) \,
\|_{\mbox{}_{\scriptstyle L^{2}(\mathbb{R}^{4})}}
^{\:\!2}
\!\;\!  +\:
C \!\!\;\!
\int_{\scriptstyle t_0}^{\;\!t}
\!\!\:\!
(\tau - t_{0})
\,
\|\,  \,\!(D^{2}\uu,D^{2}\bb) (\cdot,\tau) \,
\|_{\mbox{}_{\scriptstyle L^{2}(\mathbb{R}^{4})}}
^{\:\!2}
d\tau
\\
\leq\;
\!\!\;\!
\int_{\scriptstyle t_0}^{\;\!t}
\!\!\;\!
\,
\|\,  \,\!(D\uu,D\bb)(\cdot,\tau) \,
\|_{\mbox{}_{\scriptstyle L^{2}(\mathbb{R}^{4})}}
^{\:\!2}
d\tau,
\end{split}
\end{equation*}
for some constant $C>0$.
We proceed for general $m \geq 2 $ by induction similarly as in the case $n=3$,
\begin{multline*}
	(\:\!t - t_{0})^{m}
	\,
	\|\,  \,\!(D^{m}\uu,D^{m}\bb)(\cdot,t) \,
	\|_{\mbox{}_{\scriptstyle L^{2}(\mathbb{R}^{4})}}
	^{\:\!2}
	\\ \!\;\!+\:
	2\, \min\{\mu,\nu\} \!\!\;\!
	\int_{\scriptstyle t_0}^{\;\!t}
	\!\!\:\!
	(\tau - t_{0})^{m}
	\,
	\|\, (D^{m + 1}\uu,D^{m + 1}\bb)(\cdot,\tau) \,
	\|_{\mbox{}_{\scriptstyle L^{2}(\mathbb{R}^{4})}}
	^{\:\!2}
	d\tau
 \\
	\leq\;
	 m
	\!\!\;\!
	\int_{\scriptstyle t_0}^{\;\!t}
	\!\!\;\!
	(\tau - t_{0})^{ m -1}
	\,
	\|\,  \,\!(D^{m}\uu,D^{m}\bb)(\cdot,\tau) \,
	\|_{\mbox{}_{\scriptstyle L^{2}(\mathbb{R}^{4})}}
	^{\:\!2}
	d\tau
	\;\;\!
 \\ +
	\mbox{} \!
C
	\!\!\,\!
	\int_{\scriptstyle t_0}^{\;\!t}
	\!\!\;\!
	(\tau - t_{0})^{m}
	\,
	\|\, (D^{m + 1}\uu,D^{m + 1}\bb)(\cdot,\tau) \,
	\|_{\mbox{}_{\scriptstyle L^{2}(\mathbb{R}^{4})}}
	\!\!\;\!
	\\
	\sum_{\ell\,=\,0}^{[\,m/2\,]}
	\!\!\;\!
	\|\,  \,\! (D^{\ell}\uu,D^{\ell}\bb)(\cdot,\tau) \,
	\|_{\mbox{}_{\scriptstyle L^{4}(\mathbb{R}^{4})}}
	\,\!
	\|\,  \:\!(D^{m - \ell}\uu,D^{m - \ell}\bb)(\cdot,\tau) \,
	\|_{\mbox{}_{\scriptstyle L^{4}}}
	\!\;\!\;\!
	d\tau.
\end{multline*}
Using (2.12) and (2.5), we have
\begin{equation*}
\begin{split}
(\:\!t - t_{0})^m
\,
\|\,  \,\!(D^m\uu,D^m\bb)(\cdot,t) \,
\|_{\mbox{}_{\scriptstyle L^{2}(\mathbb{R}^{4})}}
^{\:\!2}
\!\;\!  +\:
C \!\!\;\!
\int_{\scriptstyle t_0}^{\;\!t}
\!\!\:\!
(\tau - t_{0})^m
\,
\|\,  \,\!(D^{m+1}\uu,D^{m+1}\bb) (\cdot,\tau) \,
\|_{\mbox{}_{\scriptstyle L^{2}(\mathbb{R}^{4})}}
^{\:\!2}
d\tau
\\
\leq\;
m
\!\!\;\!
\int_{\scriptstyle t_0}^{\;\!t}
\!\!\;\!
(\tau - t_0)^{m-1}
\|\,  \,\!(D^m\uu,D^m\bb)(\cdot,\tau) \,
\|_{\mbox{}_{\scriptstyle L^{2}(\mathbb{R}^{4})}}
^{\:\!2}
d\tau.
\end{split}
\end{equation*}
By the same argument in the $n=3$ case, we conclude the proof of (3.1). Now, we just have to apply a simple interpolation and the proof of Theorem I turns out.

\setcounter{section}{3}
%

%
% ************************************************
% *                                              *
% *                Section 4                     *
% *                    *                         *
% *                                              *
% ************************************************
%
\nl
% \mbox{} \vspace{-1.750cm} \\
%

{\bf 4. Proof of (1.6) } \\
\setcounter{section}{4}

We begin with $n = 2$. 
Using the Sobolev inequality (2.7$a$) 
for the pair $(\uu,\bb)$, we have
\begin{equation*}
\|\, (\uu,\bb) \,
\|_{\mbox{}_{\scriptstyle L^{\infty}(\mathbb{R}^{2})}}
\:\!\leq\;
C\|\, (\uu,\bb) \,
\|_{\mbox{}_{\scriptstyle L^{2}(\mathbb{R}^{2})}}
^{\:\!1/2}
\:\!
\|\, ( D^{2}\uu,D^{2}\bb) \,
\|_{\mbox{}_{\scriptstyle L^{2}(\mathbb{R}^{2})}}
^{\:\!1/2}.
\end{equation*}
By Main Theorem, we get
\begin{equation*}
t^{1/2}\| (\uu,\bb) \|_{L^\infty(\RR^2)} \to 0, \text{    as   } t \to \infty.
\end{equation*}
Using the same basic idea and (2.7$b$) for a pair $(\uu,\bb)$, we conclude that
\begin{equation*}
t^{n/4}\| (\uu,\bb) \|_{L^\infty(\RR^n)} \to 0, \text{    as   } t \to \infty,
\end{equation*}
for $n=2,3$.
A particular case of the fundamental Gagliardo-Nirenberg inequality ensures that
\begin{equation*}
\| (\uu,\bb) \|_{L^\infty(\RR^4)} \leq C \| (D^2\uu,D^2\bb) \|_{L^2(\RR^4)}
\end{equation*}
and using the Main Theorem again one has the same property above in dimension $n=4$. Now, we just have to apply a simple $L^2 \leftrightarrow L^\infty$ interpolation to obtain
\begin{equation*}
\begin{split}
\lim_{t \to \infty} t^{\frac{n}{4} - \frac{n}{2q}}\| (\uu,\bb)(\cdot,t) \|_{{L}^q(\RR^n)} = 0,
\end{split}
\end{equation*}
where $2 \leq q \leq \infty$ and $2 \leq n \leq 4$.
%
%
% *************************************************************
% *                                                           *
% *                      References                           *
% *                                                           *
% *************************************************************
%

%
% ----------------------------------------------------
%

\nl
\mbox{} \vspace{-0.250cm} \\
\nl
{\small

\nl
\mbox{} \vspace{-0.450cm} \\
\nl
\begin{minipage}[t]{10.00cm}
\mbox{\normalsize \textsc{Robert Guterres}} \\
Departamento de Matem\'atica Pura e Aplicada \\
Universidade Federal do Rio Grande do Sul \\
Porto Alegre, RS 91509-900, Brazil \\
E-mail: {\sf rguterres.mat@gmail.com} \\

\mbox{\normalsize \textsc{Cilon Perusato}} \\
Departamento de Matem\'atica Pura e Aplicada \\
Universidade Federal do Rio Grande do Sul \\
Porto Alegre, RS 91509-900, Brazil \\
E-mail: {\sf cilonperusato@gmail.com} \\

\mbox{\normalsize \textsc{Juliana Nunes}} \\
Departamento de Matem\'atica Pura e Aplicada \\
Universidade Federal do Rio Grande do Sul \\
Porto Alegre, RS 91509-900, Brazil \\
E-mail: {\sf juliana.s.ricardo@gmail.com} \\
\end{minipage}
}
%
% -------------------------------------------------------------
%

%
% ------------------------------------------------
%

\end{document}